\documentstyle[12pt]{article}

\begin{document}
\newcommand{\qed}{\hphantom{.}\hfill $\Box$\medbreak}
\newtheorem{exa}{Example}[section]
\newtheorem{lemma}{Lemma}[section]
\newtheorem{remark}{Remark}[section]
\newtheorem{prop}{Proposition}[section]
\newtheorem{theo}{Theorem}[section]

\def  \cS  {{\cal S}} 
\def \cG {{\cal G}}
\def \cB {{\cal B}}

\title{\bf The Doyen-Wilson theorem for $3$-sun systems\thanks{Supported
 by P.R.I.N. and I.N.D.A.M.
(G.N.S.A.G.A.)}}
\author{ {\bf Giovanni Lo Faro}\\ \small
Dipartimento di Scienze Matematiche e Informatiche, \\  \small
Scienze Fisiche e Scien\-ze della Terra \\  \small
Universit\`a di Messina, Messina, Italia\\ \small email:
lofaro@unime.it \and {\bf Antoinette Tripodi}\\ \small
Dipartimento di Scienze Matematiche e Informatiche, \\  \small
Scienze Fisiche e Scien\-ze della Terra \\  \small
Universit\`a di Messina, Messina, Italia\\ \small email:
atripodi@unime.it}
\date{}
\maketitle

\begin{abstract}
A solution to the existence problem of  $G$-designs with given subdesigns  is known  when $G$ is a triangle with $p=0,1,$ or $2$  disjoint pendent edges: for $p=0$, it is due to Doyen and Wilson, the first to pose  such a problem for  Steiner triple systems;  for $p=1$ and  $p=2$, the corresponding designs are  kite systems and bull designs, respectively. Here, a complete solution to the problem is given in  the remaining case where $G$ is a $3$-sun, i.e. a graph on six vertices consisting of  a triangle with three pendent edges which form a 1-factor.
\end{abstract}

\vspace{0.5cm}\par \noindent {\bf Keywords:} $3$-sun systems;
embedding; difference set.

  \vspace{0.5cm}\par \noindent {\bf Mathematics Subject
  Classification(2000):} 05B05, 05B30.

\section{Introduction}

If $G$ is a graph, then let $V(G)$ and $E(G)$ be the vertex-set and edge-set of $G$, respectively. The graph $K_n$ denotes the complete graph on $n$ vertices. The graph $K_m\setminus K_n$ has vertex-set $V(K_m)$ containing a distinguished subset $H$ of size $n$; the edge-set of $K_m\setminus K_n$ is $E(K_m)$ but with the $n \choose 2$ edges between the $n$ distinguished vertices of $H$ removed. This graph is sometimes referred to as a complete graph of order $m$ with a \emph{hole} of size $n$. 

Let $G$ and  $\Gamma$ be finite graphs.  A {\em $G$-design} of $\Gamma$  is a pair $(X, {\cal B})$ where $X=V(\Gamma)$  and ${\cal B}$ is a collection of isomorphic
copies of  $G$ ({\em blocks}), whose edges
partition $E(\Gamma)$. If $\Gamma=K_n$, then we refer to such a design as a \emph{$G$-design of order $n$}. 

 A $G$-design $(X_1, \cB_1)$ of order $n$ is said to be \emph{embedded} in a
$G$-design $(X_2, \cB_2)$ of order $m$ provided $X_1\subseteq X_2$ and $
\cB_1\subseteq \cB_2$ (we also say that $(X_1, \cB_1)$ is a
\emph{subdesign} (or \emph{subsystem}) of $(X_2, \cB_2)$ or $(X_2,
\cB_2)$ contains $(X_1, \cB_1)$ as subdesign). Let $N(G)$ denote the
set of integers $n$ such that there exists a $G$-design of order
$n$. A natural question to ask is: given $n,m\in N(G)$, with
$m>n$, and a $G$-design  $(X, \cB)$ of order $n$, does exists a
$G$-design of order $m$ containing $(X, \cB)$ as subdesign? Doyen
and Wilson were the first to pose this problem for $G=K_3$
(Steiner triple systems) and in 1973 they showed that given
$n,m\in N(K_3)=\{ v\equiv 1,3\!\!\!\!\pmod {6}\}$, then \emph{any
Steiner triple system of order $n$ can be embedded in a Steiner
triple system of order $m$ if and only if $m\geq 2n+1$ or $m=n$}
(see \cite{dw}). Over the years, any such problem has come to be
called a ``Doyen-Wilson problem'' and any solution a
``Doyen-Wilson type theorem''. The work along these lines is
extensive (\cite{br2}, \cite{fllhh}-\cite{hy}, \cite{lt}, \cite{lt2}, \cite{w}) and the interested
reader is referred to \cite{br} for a history of this problem.

In particular, taking into consideration the case where $G$ is a triangle with $p=0,1,2,$ or $3$ mutually disjoint pendent edges,  a solution to the Doyen-Wilson problem is known  when $p=0$  (Steiner triple systems, \cite{dw}),  $p=1$  (kite systems, \cite{lt,lt2}) and  $p=2$  (bull designs, \cite{fllhh}). Here, we deal with the remaining case ($p=3$) 
where  $G$ is a  $3$-\emph{sun}, i.e. a graph on six vertices consisting of a triangle with three pendent edges which form a 1-factor, by giving a complete solution to the Doyen-Wilson problem for $G$-designs where $G$ is a $3$-sun (\emph{$3$-sun systems}).

\section{Notation and basic lemmas}

The  \emph{$3$-sun} consisting of the triangle $(a, b, c)$ and the  three  disjoint pendent edges $\{{a, d\}, \{b, e\}, \{c, f\}}$  is denoted   by  $(a,b,c;d,e,f)$. A  $3$-sun system of order $n$ (briefly, 3SS$(n)$) exsits if and only if $n\equiv 0,1,4,9\!\!\!\!\pmod {12}$  and if $(X, \cS)$ is a 3SS$(n)$, then $|\cS|=\frac {n(n-1)}{12}$ (see \cite{yg}).

 Let $n,m \equiv 0,1,4,9\!\!\!\!\pmod {12}$, with $m=u+n$,  $u\geq 0$. The Doyen-Wilson problem for $3$-sun systems  is equivalent to the existence problem of decompositions of $K_{u+n}\setminus K_n$ into $3$-suns. 

Let $r$ and $s$ be  integers with $r<s$, define $[r, s] = \{r, r + 1, . . . , s\}$ and $[s, r] = \emptyset$. Let $Z_u=[0, u-1] $ and  $H=\{\infty _1, \infty _2,\ldots ,\infty _t\}$, $H\cap Z_u=\emptyset$. If $S=(a,b,c; d,e,f)$ is a 3-sun whose vertices belong to  $Z_u \cup H$ and $i\in Z_u$,    let $S+i=(a+i,b+i,c+i; d+i,e+i,f+i)$, where the sums are modulo $u$ and $\infty +i=\infty$, for every $\infty \in H$.  The set $(S)= \{S+i : i\in Z_u\}$  is called the
\emph{orbit of $S$ under $Z_u$} and $S$ is  a \emph{base block} of $(S)$. 

To solve the Doyen-Wilson problem for $3$-sun systems   we  use the \emph{difference method} (see  \cite{ls}, \cite{p}). For every pair of distinct elements  $i,j \in Z_u$,  define $|i-j|_u $= min$\{|i-j|, u-|i-j|\}$ and set $D_u =\{|i-j|_u : i,j\in Z_u\}=\{1, 2, \ldots, \lfloor \frac{u}{2}\rfloor\}$. The elements of $D_u$ are called \emph{differences} of $Z_u$. 
 For
any $d\in D_u$, $d\neq \frac u2$,  we can form a single
2-factor $\{\{i, d+i\}\, :\, i\in Z_u\}$, while if $u$ is even and  $d=
\frac u2$, then we can form a 1-factor $\{\{i, i+\frac u2\}\, :\,
0 \leq i\leq \frac u2 -1 \}$. It is also worth remarking that
2-factors obtained from distinct differences are disjoint from
each other and from the 1-factor.

If  $D\subseteq D_u$, denote by $\langle Z_u\cup H, D\rangle$
the graph with vertex-set  $V=Z_u\cup H$ and the edge-set $
E=\{\{i, j\}\, :\,|i-j|_u=d ,\ d\in
D \}\cup \{\{\infty , i\}\, :\,  \infty \in H,\ i\in Z_u\}$. The graph $\langle Z_u\cup H, D_u\rangle$ is the complete graph $K_{u+t}\setminus K_t$ based on $Z_u\cup H $ and  having $H$ as hole. The elements of $H$ are called  \emph{infinity points}.

Let $X$ be a set of size  $n \equiv 0,1,4,9\!\!\!\!\pmod {12}$. The aim of the paper is to decompose the graph $\langle Z_u\cup X, D_u\rangle$ into $3$-suns. To obtain our main result the  $\langle Z_u\cup X, D_u\rangle$ will be regarded as a union of  suitable edge-disjoint subgraphs of type $\langle Z_u\cup H, D\rangle$ (where $H\subseteq X$ may be empty, while  $D\subseteq D_u$ is always non empty) and  then each subgraph will be decomposed  into 3-suns  by using the  lemmas given in this section.  
From here on suppose $u\equiv 0,1,3,4,5,7,8,$ $9,11\!\!\!\!\pmod {12}$.

\bigskip 

Lemmas \ref{2-1} - \ref{8-2} give  decompositions of subgraphs of type $\langle Z_u\cup H, D\rangle$  where $D$ contains particular differences, more precisely, $D=\{ 2\}$, $D=\{ 2,4\}$ or $D=\{ 1, $ $\frac u3\}$.

\begin{lemma} \label{2-1}
Let $u\equiv 0 \pmod {4}$, $u\geq 8$. Then the  graph $\langle Z_u\cup \{\infty _1, \infty _2\}, $ $\{ 2\}\rangle$
 can be decomposed into $3$-suns.
\end{lemma}
{\em Proof.} Consider the $3$-suns 
$$\begin{array}{l}
 (\infty _1,{2+4i},{4i}; {3+4i},4+4i,\infty _2), \\
 (\infty _2,{3+4i},{1+4i}; {2+4i},5+4i,\infty _1) , 
\end{array}$$ for $i=0,1, \ldots, \frac{u}{4}-1$. \hfill$\Box$
\bigskip

\begin{lemma} \label{4-1}
Let $u\equiv 0 \pmod {12}$. Then the graph $\langle Z_u\cup \{\infty _1, \infty _2, \infty _3, $ $\infty _4\}, \{ 2\}\rangle$
 can be decomposed into $3$-suns.
\end{lemma}
{\em Proof.} Consider the $3$-suns 
$$\begin{array}{l}
 (\infty _1,{12i},{2+12i}; {7+12i},,\infty _3,\infty _4), \\
 (\infty _1,{4+12i},{6+12i}; {9+12i}, \infty _3,\infty _4), \\
 (\infty _1,{8+12i},{10+12i}; {11+12i}, \infty _3,\infty _4), \\
 (\infty _2,{2+12i},{4+12i}; {1+12i}, \infty _3,\infty _4), \\
 (\infty _2,{6+12i},{8+12i}; {7+12i}, \infty _3,\infty _4), \\
 (\infty _2,{10+12i},{12+12i}; {11+12i}, \infty _3,\infty _4),\\
 (\infty _3,{1+12i},{3+12i}; {9+12i}, \infty _1,\infty _2),\\
 (\infty _3,{5+12i},{7+12i}; {11+12i}, \infty _1,9+12i),\\
 (\infty _4,{3+12i},{5+12i}; {1+12i}, \infty _1,\infty _2),\\
 (\infty _4,{9+12i},{11+12i}; {7+12i}, \infty _2,13+12i),\end{array}$$ for $i=0,1, \ldots, \frac{u}{12}-1$. \hfill$\Box$
\bigskip

\begin{lemma} \label{4-2}
The graph $\langle Z_u\cup \{\infty _1, \infty _2,\infty _3, $ $\infty _4\}, \{2,4\}\rangle$,  $u\geq 7$, $u\neq 8$, can be decomposed into $3$-suns.
\end{lemma}

{\em Proof.} Let $u=4k+r$, with $r=0,1,3$, and consider the $3$-suns
$$
\begin{array}{l}
(\infty _1,{4+4i},{6+4i}; {5+4i},8+4i,\infty _4),\\ 
(\infty _2,{5+4i},{7+4i}; {6+4i},9+4i,\infty _1),\\ 
(\infty _3,{6+4i},{8+4i}; {7+4i},10+4i,\infty _2),\\ 
(\infty _4,{7+4i},{9+4i}; {8+4i},11+4i,\infty _3),
\end{array}$$ 
\noindent for $i=0,1, \ldots, k-3$,  $k\geq3$, plus the following blocks as the case may be.\\
If $r=0$,
$$\begin{array}{l}
(\infty _1,{0},{2}; {1},4,\infty _4),\\ 
(\infty _2,{1},{3}; {2},5,\infty _1),\\ 
(\infty _3,{2},{4}; {3},6,\infty _2),\\ 
(\infty _4,{3},{5}; {4},7,\infty _3),\\
(\infty _1,{4k-4},{4k-2}; {4k-3},0,\infty _4),\\ 
(\infty _2,{4k-3},{4k-1}; {4k-2},1,\infty _1),\\ 
(\infty _3,{4k-2},{0}; {4k-1},2,\infty _2),\\ 
(\infty _4,{4k-1},{1}; {0},3,\infty _3).
\end{array}$$ 
If $r=1$,
$$\begin{array}{l}
(\infty _1,{0},{2}; {1},4,\infty _2),\\ 
(\infty _2,{1},{3}; {0},5,\infty _1),\\ 
(\infty _3,{2},{4}; {3},6,\infty _2),\\ 
(\infty _4,{3},{5}; {4},7,\infty _3),\\
(\infty _1,{4k-4},{4k-2}; {4k-3},4k,\infty _2),\\ 
(\infty _2,{4k-3},{4k-1}; {4k},0,\infty _1),\\ 
(\infty _3,{4k-2},{4k}; {4k-1},1,\infty _1),\\ 
(\infty _4,{4k-1},{0}; {4k-2},2,\infty _3),\\ 
(\infty _4,{4k},{1}; {2},3,\infty _3).
\end{array}$$ 
If $r=3$,
$$\begin{array}{l}
(\infty _1,{0},{2}; {1},4,\infty _4),\\ 
(\infty _2,{1},{3}; {2},5,\infty _1),\\ 
(\infty _3,{2},{4}; {3},6,\infty _2),\\ 
(\infty _4,{3},{5}; {4},7,\infty _3),\\
(\infty _1,{4k-4},{4k-2}; {4k-3},4k,\infty _4),\\ 
(\infty _2,{4k-3},{4k-1}; {4k-2},4k+1,\infty _1),\\ 
(\infty _3,{4k-2},{4k}; {4k-1},4k+2,\infty _2),\\ 
(\infty _4,{4k-1},{4k+1}; {4k},0,\infty _3),\\ 
(\infty _1,{4k},{4k+2}; {4k+1},1,\infty _4),\\ 
(\infty _2,{4k+1},{0}; {4k+2},2,\infty _4),\\ 
(\infty _3,{4k+2},{1}; {0},3,\infty _4). 
\end{array}$$ 
With regard to the difference $4$ in $Z_7$, note that $|4|_7=3$ and the seven distinct blocks  obtained for $k=1$  and $r=3$ 
gives a decomposition of $\langle Z_7\cup \{\infty _1, \infty _2,\infty _3, $ $\infty _4\}, \{2,3\}\rangle$ into $3$-suns.
\hfill$\Box$

\begin{lemma} \label{8-2}
Let $u\equiv 0 \pmod {3}$, $u\geq 12$. Then the graph $\langle Z_u\cup \{\infty _1, \infty _2,$ $\ldots, \infty _8\}, \{ 1, $ $\frac u3\}\rangle$
 can be decomposed into $3$-suns.
\end{lemma}
{\em Proof.}  
If $u\equiv 0 \pmod {6}$ consider the $3$-suns:
$$\begin{array}{l}
(\infty _1,{2i},{\frac u3+2i}; {2\frac u3+2i},\infty _5,\infty _6), \ i=0,1, \ldots, \frac u6-1, \\
(\infty _1,{1+2i},{\frac u3+1+2i}; {2\frac u3+1+2i},\infty _6,\infty _5), i=0,1,\  \ldots, \frac u6-1, \\
(\infty _2,{2\frac u3+2i},{\frac u3+2i}; {2+2i},2i,\infty _5), \ i=0,1, \ldots, \frac u6-2, \\
(\infty _2,{2\frac u3+1+2i},{\frac u3+1+2i}; {3+2i},1+2i,\infty _6), i=0,1,\  \ldots, \frac u6-2, \\
(\infty _2,{u-2},{2\frac u3-2}; {0},\frac u3-2,\infty _5), \  
(\infty _2,{u-1},{2\frac u3-1}; {1},\frac u3-1,\infty _6), \\
(\infty _3,{2i},{1+2i}; {2\frac u3+2i},\infty _7,\infty _8), \ i=0,1, \ldots, \frac u6-1, \\
(\infty _3,{\frac u3+2i},{\frac u3+1+2i}; {2\frac u3+1+2i},\infty _7,\infty _8), i=0,1,\  \ldots, \frac u6-1, \\
(\infty _4,{1+2i},{2+2i}; {2\frac u3+2+2i},\infty _7,\infty _8), \ i=0,1, \ldots, \frac u6-1, \\
(\infty _4,{\frac u3+1+2i},{\frac u3+2+2i}; {2\frac u3+1+2i},\infty _7,\infty _8), i=0,1,\  \ldots, \frac u6-1, \\
(\infty _5,{2\frac u3+2i},{2\frac u3+1+2i}; {1+2i},\infty _7,\infty _8), \ i=0,1, \ldots, \frac u6-1, \\
(\infty _6,{2\frac u3+3+2i},{2\frac u3+4+2i}; {2+2i},\infty _7,\infty _8), i=0,1,\  \ldots, \frac u6-2, \\
(\infty _6,{2\frac u3+1},{2\frac u3+2}; {2\frac u3},\infty _7,\infty _8).
\end{array}$$ 
\noindent 
If $u\equiv 3 \pmod {6}$ consider the $3$-suns:
$$\begin{array}{l}
(\infty _1,{2i},{\frac u3+2i}; {2\frac u3+2i},\infty _5,\infty _6), \ i=0,1, \ldots, \frac{u-3}{6}, \\
(\infty _1,{1+2i},{\frac u3+1+2i}; {2\frac u3+1+2i},\infty _6,\infty _5), i=0,1,\  \ldots, \frac{u-9}{6}, \\
(\infty _2,{2\frac u3+2i},{\frac u3+2i}; {2+2i},2i,\infty _5), \ i=0,1, \ldots, \frac{u-9}{6}, \\
(\infty _2,{u-1},{2\frac u3-1}; {0},\frac u3-1,\infty _5), \\
(\infty _2,{2\frac u3+1+2i},{\frac u3+1+2i}; {3+2i},1+2i,\infty _6), i=0,1,\  \ldots, \frac{u-15}{6}, \\
(\infty _2,{u-2},{2\frac u3-2}; {1},\frac u3-2,\infty _6), \\  
(\infty _3,{2i},{1+2i}; {2\frac u3+2i},\infty _7,\infty _8), \ i= 2,3,\ldots, \frac{u-3}{6}, \\
(\infty _3,{0},{1}; {2\frac u3},  \infty _6,\infty _8),  
(\infty _3,{2},{3}; {2\frac u3+2},  \infty _6,\infty _8), \\  
(\infty _3,{\frac u3+1+2i},{\frac u3+2+2i}; {2\frac u3+1+2i},\infty _7,\infty _8), i=0,1,\  \ldots, \frac{u-9}{6}, \\
(\infty _4,{1+2i},{2+2i}; {2\frac u3+2+2i},\infty _7,\infty _8), \ i=0,1, \ldots, \frac{u-9}{6}, \\
(\infty _4,{\frac u3+2i},{\frac u3+1+2i}; {2\frac u3+1+2i},\infty _7,\infty _8), i=0,1,\  \ldots, \frac{u-3}{6}, \\
(\infty _5,{2\frac u3+2i},{2\frac u3+1+2i}; {1+2i},\infty _7,\infty _8), \ i=0,1, \ldots, \frac{u-9}{6}, \\
(\infty _6,{2\frac u3+1+2i},{2\frac u3+2+2i}; {4+2i},\infty _7,\infty _8), \ i=0,1,\  \ldots, \frac{u-15}{6},\\
(\infty _6,{u-2},{u-1}; {2\frac u3},\infty _7,\infty _8), \
(\infty _7,{u-1},{0}; {2},\infty _5,\infty _8).
\end{array}$$ 
\hfill$\Box$
\bigskip

Lemmas \ref{3-2u2} - \ref{3-3u2} allow  to  decompose  $\langle Z_u\cup H, D\rangle$ where $u$ is even and $ D$ contains the difference $\frac u2$.

\begin{lemma} \label{3-2u2}
Let $u$ be even, $u\geq 8$. Then the graph $\langle Z_u\cup \{\infty _1, \infty _2, \infty _3\}, \{1,$ $\frac u2\}\rangle$
 can be decomposed into $3$-suns.
\end{lemma}
{\em Proof.} Consider the $3$-suns 
$$\begin{array}{l}
 (\infty _1,{2i},{1+2i}; {\frac u2+2+2i},\frac u2+2i,\infty _3), \ i=0,1, \ldots, \frac{u}{4}-2, \\
(\infty _1,{\frac u2-2},{\frac u2-1}; {\frac u2},u-2,\infty _3),\\
(\infty _2,{1+2i},{\frac u2+1+2i}; {2i},2+2i,\infty _1), \ i=0,1, \ldots, \frac{u}{4}-1, \\
(\infty _3,{\frac u2+1+2i},{\frac u2+2i}; {2i},\frac u2+2+2i,\infty _2), \ i=0,1, \ldots, \frac{u}{4}-1.
\end{array}$$ 
\hfill$\Box$
\bigskip

\begin{lemma} \label{4-2u2}
Let $u\equiv 0 \pmod {12}$. Then the graph $\langle Z_u\cup \{\infty _1, \infty _2, \infty _3, \infty _4\}, $ $\{1,\frac u2\}\rangle$
 can be decomposed into $3$-suns.
\end{lemma}
{\em Proof.} Consider the $3$-suns 
$$\begin{array}{l}
 (\infty _1,{6i},{\frac u2+6i}; {4+6i},\infty _3 ,\infty _2), \\  
 (\infty _1,{1+6i},{\frac u2+1+6i}; {5+6i},\infty _4 ,\infty _2),\\
 (\infty _1,{2+6i},{\frac u2+2+6i}; {\frac u2+3+6i},\infty _4,\infty _3), \\  
 (\infty _2,{1+6i},{ 6i}; {\frac u2+3+6i},\infty _3 ,\infty _4),\\
 (\infty _2,{2+6i},{3+6i}; {\frac u2+4+6i},1+6i,\infty _4), \\  
 (\infty _2,{5+6i},{ 4+6i}; {\frac u2+5+6i},6+6i ,3+6i), \\ (\infty _3,{3+6i},{\frac u2+3+6i}; {2+6i},\infty _1 ,\frac u2+2+6i), \\  
 (\infty _3,{4+6i},{\frac u2+4+6i}; {\frac u2+6i},\infty _4 ,\infty _1),\\
  (\infty _3,{5+6i},{\frac u2+5+6i}; {\frac u2+1+6i},\infty _4 , \infty _1), \\  
 (\infty _4,{\frac u2+1+6i},{\frac u2+2+6i}; {\frac u2+3+6i},\frac u2+6i ,\infty _2), \\
  (\infty _4,{\frac u2+4+6i},{\frac u2+5+6i}; {\frac u2+6i},\frac u2+3+6i ,\frac u2+6+6i), \end{array}$$ for $i=0,1, \ldots, \frac{u}{12}-1$. \hfill$\Box$
\bigskip

\begin{lemma} \label{6-2u2}
Let $u$ be even, $u\geq 8$. Then the graph $\langle Z_u\cup \{\infty _1, \infty _2, \dots, \infty _6\}, \{1,\frac u2\}\rangle$
 can be decomposed into $3$-suns.
\end{lemma}
{\em Proof.} Consider the $3$-suns 
$$\begin{array}{l}
 (\infty _1,{2i},{1+2i}; {\frac u2+2+2i},\frac u2+2i,\infty _3), i=0,1, \ldots, \frac{u}{4}-2,\\
 (\infty _1,{\frac u2-2},{\frac u2-1}; {\frac u2},u-2,\infty _3), \\
 (\infty _2,{1+2i},{\frac u2+1+2i}; {2i},\infty _6,\infty _1), i=0,1, \ldots, \frac{u}{4}-1,\\
 (\infty _3,{\frac u2+1+2i},{\frac u2+2i}; {2i},\infty _6,\infty _2),  i=0,1, \ldots, \frac{u}{4}-1,\\
 (\infty _4,{1+2i},{2+2i}; {\frac u2+2+2i},\infty _5,\infty _6),  i=0,1, \ldots, \frac{u}{4}-1,\\
 (\infty _5,{\frac u2+1+2i},{\frac u2+2+2i}; {2+2i},\infty _4,\infty _6),   i=0,1, \ldots, \frac{u}{4}-1.
\end{array}$$ \hfill$\Box$
\bigskip

\begin{lemma} \label{7-2u2}
Let  $u\equiv 0 \pmod {12}$. Then the graph $\langle Z_u\cup \{\infty _1, $ $\infty _2, \dots, \infty _7\}, \{1,\frac u2\}\rangle$
 can be decomposed into $3$-suns.
\end{lemma}
{\em Proof.} Consider the $3$-suns 
$$\begin{array}{l}
 (\infty _1,{6i},{\frac u2+6i}; {4+6i},\infty _7 ,\infty _2), \\  
 (\infty _1,{1+6i},{\frac u2+1+6i}; {\frac u2+3+6i},\infty _7 ,\infty _4),\\
 (\infty _1,{2+6i},{\frac u2+2+6i}; {\frac u2+5+6i},\infty _5,\infty _2), \\  
 (\infty _2,{3+6i},{ \frac u2+3+6i}; {6i},\infty _1 ,\infty _4),\\
 (\infty _2,{4+6i},{\frac u2+4+6i}; {2+6i}, \infty _7,\infty _1), \\  
 (\infty _2,{5+6i},{ \frac u2+5+6i}; {\frac u2+1+6i},  \infty _1,\infty _7),  \\ 
 (\infty _3,{ 6i},{1+6i}; {\frac u2+6i},\infty _5 ,\infty _6), \\  
 (\infty _3,{2+6i},{3+6i}; {\frac u2+2+6i},\infty _7 ,\infty _6),\\
  (\infty _3,{4+6i},{ 5+6i}; {\frac u2+5+6i},\infty _5 , \infty _6), \\  
 (\infty _4,{ 1+6i},{\ 2+6i}; {\frac u2+6+6i},\infty _2,\infty _6), \\
  (\infty _4,{3+6i},{4+6i}; {\frac u2+4+6i},\infty _7 ,\infty _6),\\
  (\infty _4,{5+6i},{6+6i}; {\frac u2+5+6i},\infty _7 ,\infty _6),\\  
 (\infty _5,{\frac u2+ 6i},{ \frac u2+1+6i}; { 1+6i},  \infty _7,\infty _3), \\  
 (\infty _5,{\frac u2+2+6i},{ \frac u2+3+6i}; { 3+6i},  \infty _7,\infty _3),  \\  
 (\infty _5, {\frac u2+4+6i},{ \frac u2+5+6i}; { 5+6i},  \infty _7,\frac u2+6+6i),\\  
 (\infty _6,{\frac u2+1+6i},{ \frac u2+2+6i}; { \frac u2+5+6i},  \infty _7,\infty _4), \\  
 (\infty _6,{\frac u2+3+6i},{ \frac u2+4+6i}; { \frac u2+6+6i},  \infty _7,\infty _3), 
 \end{array}$$ for $i=0,1, \ldots, \frac{u}{12}-1$. \hfill$\Box$
\bigskip

\begin{lemma} \label{3-3u2}
Let  $u\equiv 0 \pmod {12}$. Then the graph $\langle Z_u\cup \{\infty _1, \infty _2, \infty _3\}, $ $\{1,2,\frac u2\}\rangle$
 can be decomposed into $3$-suns.
\end{lemma}
{\em Proof.} Consider the $3$-suns 
$$\begin{array}{l}
 (\infty _1,{6i},{1+6i}; {\frac u2+1+6i},\frac u2+6i ,3+6i), \\  
 (\infty _1,{2+6i},{3+6i}; {\frac u2+5+6i},\frac u2+2+6i ,5+6i), \\  
 (\infty _1,{4+6i},{5+6i}; {\frac u2+2+6i},\frac u2+4+6i ,7+6i), \\
  (\infty _1,{\frac u2+3+6i},{\frac u2+4+6i}; {\frac u2+6i},\frac u2+2+6i , \infty _2), \\  
 (\infty _2,{1+6i},{ \frac u2+1+6i}; {\frac u2+3+6i}, 2+6i, \frac u2+2+6i),\\
 (\infty _2,{3+6i},{4+6i}; {2+6i},\frac u2+3+6i, 6+6i), \\  
 (\infty _2,{5+6i},{ \frac u2+5+6i}; {\frac u2+2+6i},6+6i ,\frac u2+6+6i), \\
 (\infty _3,{2+6i},{  6i}; {1+6i}, 4+6i, \infty _2 ), \\  
 (\infty _3,{\frac u2+2+6i},{\frac u2 +6i}; {4+6i},\frac u2+4+6i ,\infty _2),\\
 (\infty _3,{\frac u2+1+6i},{\frac u2+3+6i}; {3+6i},\frac u2+6i, \frac u2+5+6i), \\  
 (\infty _3,{\frac u2+5+6i},{\frac u2+4+6i}; {5+6i},\frac u2+7+6i, \frac u2+6+6i), \end{array}$$ 
  for $i=0,1, \ldots, \frac{u}{12}-1$. \hfill$\Box$
\bigskip

The following lemma ``\,combines\,''  one infinity point with one difference $d\neq \frac u2, \frac u3$ such that $\frac u{gcd(u,d)}\equiv 0 \pmod {3}$ (therefore, $u\equiv 0 \pmod {3}$).

\begin{lemma} \label{1-1}
 Let $u\equiv 0 \pmod {3}$ and $d\in D_u\setminus \{\frac u2, \frac u3\}$  such that 
 $p=\frac u{gcd(u,d)} \equiv 0 \pmod {3}$. Then the graph $\langle Z_u\cup \{\infty\}, \{d\}\rangle$
 can be decomposed into $3$-suns.
\end{lemma}
{\em Proof.} The subgraph $\langle Z_u, \{d\}\rangle$  can be decomposed into $\frac{u}{p}$ cycles of length $p=3q$, $q \geq 2$. \\
If $q>2$, let
$(x_1,x_2, \ldots, x_{3q})$ be a such cycle and consider the $3$-suns
$$
(\infty ,x_{2+3i},x_{3+3i}; x_{7+3i},x_{1+3i},x_{4+3i}),$$
\noindent for $i=0,1, \ldots, q-1$ (where the sum is modulo $3q$). \\
If $q=2$,  let $(x_1^{(j)},x_2^{(j)},x_3^{(j)},x_4^{(j)},x_5^{(j)} , x_{6}^{(j)})$, $j=0,1,\ldots,\frac{u}{6}-1$, be the $6$-cycles decomposing  $\langle Z_u, \{d\}\rangle$ and  consider the $3$-suns
$$ (\infty ,x_{2}^{(j)},x_{3}^{(j)}; x_{1}^{( j+1)},x_{1}^{(j)},x_{4}^{(j)}),\  (\infty ,x_{5}^{(j)},x_{6}^{(j)}; x_{4}^{(j+1)},x_{4}^{(j)},x_{1}^{(j)}),$$
\noindent for $j=0,1, \ldots, \frac{u}{6}-1$ (where the sums are modulo $\frac{u}{6}$). \hfill$\Box$
\bigskip

Subsequent Lemmas \ref{5-1} - \ref{1-5} allow to decompose   $\langle Z_u\cup H, D\rangle$, where $|H|=1,2,3,5$, $|D|=6-|H|$ and $\frac u2 \not \in D$; here, $u$ and $D$ are any with the unique condition that if $D$ contains at least three differences  $d_1,d_2,d_3$, then $d_3=d_2-d_1 $  or $d_1+d_2+ d_3=u$.

\begin{lemma} \label{1-5}
Let $d_1,d_2,d_3,d_4,d_5 \in D_u\setminus \{\frac u2\}$  such that
$d_3=d_2-d_1 $ or $d_1+d_2+ d_3=u$. Then  the graph $\langle Z_u\cup \{\infty \},
\{d_1,d_2,d_3,d_4,d_5\}\rangle$ can be decomposed into $3$-suns.
\end{lemma}

{\em Proof.} If $d_3=d_2-d_1 $, consider the orbit  of $(d_1,d_2, 0;\infty , d_2+ d_5  ,d_4) $ 
(or $(d_1,d_2, 0;\infty , d_2+ d_5 ,-d_4) $, if 
$d_2+d_5=d_4$) under $Z_u$. If $d_1+d_2+ d_3=u$, consider the orbit  of $(-d_1,d_2, 0;\infty , d_2+ d_5  ,d_4) $ (or $(-d_1,d_2, 0;\infty , d_2+ d_5 ,-d_4) $, if 
$d_2+d_5=d_4$) under $Z_u$. \hfill$\Box$
\bigskip

\begin{lemma} \label{2-4}
Let $d_1,d_2,d_3,d_4 \in D_u\setminus \{\frac u2\}$  such that
$d_3=d_2-d_1 $ or $d_1+d_2+ d_3=u$. Then  the graph $\langle Z_u\cup \{\infty _1, \infty _2\},
\{d_1,d_2,d_3,d_4\}\rangle$ can be decomposed into $3$-suns.
\end{lemma}

{\em Proof.} Consider the orbit  of  $(d_1,d_2, 0;\infty_1 , \infty_2 ,d_4) $ or  $(-d_1,d_2, 0;\infty_1 , $ $\infty_2 ,d_4) $ under $Z_u$ when, respectively,  $d_3=d_2-d_1 $ or $d_1+d_2+ d_3=u$. \hfill$\Box$
\bigskip

\begin{lemma} \label{3-3}
Let $d_1,d_2,d_3\in D_u\setminus \{\frac u2\}$  such that
$d_3=d_2-d_1 $ or $d_1+d_2+ d_3=u$. Then  the graph $\langle Z_u\cup \{\infty _1, \infty _2,\infty _3\},
\{d_1,d_2,d_3\}\rangle$ can be decomposed into $3$-suns.
\end{lemma}

{\em Proof.} Consider the orbit  of $(d_1,d_2, 0;\infty_1 , \infty_2 ,\infty_3) $  or $(-d_1,d_2, 0;\infty_1 , $ $\infty_2 ,\infty_3) $  under $Z_u$ when, respectively, $d_3=d_2-d_1 $ or $d_1+d_2+ d_3=u$.\hfill $\Box$
\bigskip

\begin{lemma} \label{5-1}
Let  $d\in D_u\setminus \{\frac u2\}$, the graph $\langle Z_u\cup
\{\infty _1, \infty _2,$ $\infty _3, \infty _4, \infty _5\}, \{d\}\rangle$
 can be decomposed into $3$-suns.
\end{lemma}

{\em Proof.} The subgraph $\langle Z_u, \{d\}\rangle$ is regular
of degree 2 and so can be decomposed into $l$-cycles, $l\geq 3$. Let
$(x_1,x_2, \ldots, x_{l})$ be a such cycle. Put $l=3q+r$, with $r=0,1,2$, and consider the $3$-suns with the sums  modulo $l$
$$
\begin{array}{l}
(\infty _1,x_{1+3i},x_{2+3i}; x_{3+3i},\infty _4,\infty _5),\\ (\infty _2,x_{2+3i},x_{3+3i}; x_{4+3i},\infty _4,\infty _5), \\ (\infty _3,x_{3+3i},x_{4+3i}; x_{5+3i},\infty _4,\infty _5), 
\end{array}$$ 
\noindent for $i=0,1, \ldots, q-2$,  $q \geq 2$, plus the following blocks as the case may be.\\
If $r=0$,
$$\begin{array}{l}
(\infty _1,x_{3q-2},x_{3q-1}; x_{3q},\infty _4,\infty _5),\\
(\infty _2,x_{3q-1},x_{3q}; x_{1},\infty _4,\infty _5),\\
(\infty _3,x_{3q},x_{1}; x_{2},\infty _4,\infty _5).
\end{array}$$ 
If $r=1$,
$$\begin{array}{l}
(\infty _1,x_{3q-2},x_{3q-1}; x_{3q+1},\infty _4,\infty _5),\\
(\infty _2,x_{3q-1},x_{3q}; x_{1},\infty _4,\infty _1),\\
(\infty _3,x_{3q},x_{3q+1}; x_{2},\infty _4,\infty _2),\\
(\infty _5,x_{3q+1},x_{1}; x_{3q},\infty _4,\infty _3).
\end{array}$$ 
If $r=2$,
$$\begin{array}{l}
(\infty _1,x_{3q-2},x_{3q-1}; x_{3q+2},\infty _4,\infty _5),\\
(\infty _2,x_{3q-1},x_{3q}; x_{1},\infty _4,\infty _5),\\
(\infty _3,x_{3q},x_{3q+1}; x_{2},\infty _1,\infty _2),\\
(\infty _4,x_{3q+1},x_{3q+2}; x_{3q},\infty _1,\infty _3),\\
(\infty _5,x_{3q+2},x_{1}; x_{3q+1},\infty _2,\infty _3).
\end{array}$$  \hfill $\Box$
\bigskip

Finally, after settling the infinity points by using the above lemmas, if $u$ is large  we need to decompose the  subgraph  $ \langle Z_u, L \rangle$, where $L$ is the set of the differences unused  (\emph{difference leave}). Since by applying Lemmas \ref{2-1}-\ref{3-3} it could be necessary to use the differences $1$, $2$ or $4$, while Lemma \ref{5-1} does not impose any restriction,  it is possible to combine  infinity points and differences in such a way that the difference leave $L$ is the set of the ``\,small\,''  differences, where  $1$, $2$ or $4$ could possibly be avoided.

\begin{lemma} \label{55}
Let $\alpha \in \{0,4,8\}$ and  $u$, $s$ be positive integers such that $u> 12s+\alpha$. Then there exists a decomposition of $ \langle Z_u, L \rangle$  into   $3$-suns, where:
\begin{itemize} 
\item[$i)$]  $\alpha=0$ and $L= [1,6s]$; 
\item[$ii)$] $\alpha=4$ and $L= [3,6s+2]$; 
\item[$iii$)]  $\alpha=8$ and  $L= [3,6s+4]\setminus \{4,6s+3\}$.
 \end{itemize}
\end{lemma}

{\em Proof.} 
\begin{itemize} 
\item[$i)$]  Consider the orbits $(S_j)$ under  $Z_u$, where
$ S_j=(5s+1+j,\ 5s-j,\ 0; 3s,\ s, u-2-2j)$, $j=0,1,\ldots,s-1$.
\item[$ii)$] Consider the orbits in $i)$, where $(S_0)$ is replaced  with the orbit of $(6s+1,\ 4s,\ 0; s,$ $\ 9s, 6s+2)$.
\item[$iii$)]  Consider the orbits in $i)$, where the orbits $(S_0)$ and $(S_1)$ are replaced with the orbits of $(6s+1,\ 4s,\ 0; s,\ 9s, 6s+4)$ and $(5s+2,\ 5s-1,\ 0; 3s,\ s, 6s+2)$.\hfill$\Box$
\end{itemize}

\bigskip


\section{The main result}

Let $(X, \cS)$ be a  $3$-sun system  of order $n$ and $m\equiv
0,1,4,9 \pmod {12}$.

\begin{lemma} \label{NC}
If  $(X, \cS)$  is embedded in a  $3$-sun system  of order
$m>n$, then $m\geq \frac 75 n+1$.
\end{lemma}

{\em Proof.} Suppose $(X, \cS)$ embedded in $(X', \cS')$, with
$|X'|=m=n+u$ ($u$ positive integer). Let $c_i$ be the number of 3-suns of $\cS'$ each of which contains
exactly $i$ edges in $X'\setminus X$. Then $\sum_{i=1}^6i\times c_i=
{u\choose 2 }$ and $\sum_{i=1}^5(6-i)c_i= u\times n$, from which it follows
$6c_2+12c_3+18c_4+24c_5+30c_6= \frac {u(5u-2n-5)}2$ and so $u\geq \frac 25
n+1$ and  $m\geq \frac 75 n+1$. \hfill$\Box$

\bigskip

By previous Lemma:
\begin{enumerate}
\item if $n= 60k+5r$, $r=0,5,8,9$, then $m\geq 84k+7r+1$;
\item if $n= 60k+5r+1$, $r=0,3,4,7$, then $m\geq  84k+7r+3$;
\item if $n= 60k+5r+2$, $r=2,7,10,11$, then $m\geq  84k+7r+4$;
\item if $n= 60k+5r+3$, $r=2,5,6,9$, then $m\geq 84k+7r+6$;
\item if $n= 60k+5r+4$, $r=0,1,4,9$, then $m\geq 84k+7r+7$.
\end{enumerate}
\bigskip

In order to prove that the necessary conditions for embedding a  $3$-sun system $(X, \cS)$ of order $n$ in a $3$-sun  system of order $m=n+u$, $u>0$ are also sufficient, the graph $\langle Z_u\cup X, D_u \rangle$ will be expressed as a union of edge-disjoint subgraphs  $\langle Z_u\cup X, D_u \rangle=\langle Z_u\cup X, D \rangle  \cup  \langle Z_u, L \rangle$, where $L=D_u\setminus D$ is the difference leave, and $\langle Z_u\cup X, D \rangle$ (if necessary, expressed itself as a union of subgrapphs)  will be decomposed by using Lemmas \ref{2-1}-\ref{5-1}, while if $L\neq \emptyset$,  $ \langle Z_u, L \rangle$ will be decomposed by Lemma \ref{55}. To obtain our main result we will distiguish the five cases $1.-5.$ listed before by giving a general proof  for any  $k\geq 0$ with the exception of a few cases for $k=0$, which  will be indicated by a star $\star$ and solved in Appendix.  Finally, note that:
\begin{itemize}
\item[a)] $u\equiv 0,1,4,$ or $9 \pmod {12}$, if $n\equiv 0 \pmod {12}$;
\item[b)] $u\equiv 0,3,8,$ or $11\pmod {12}$, if $n\equiv 1 \pmod {12}$;
\item[c)] $u\equiv 0,5,8,$ or $9\pmod {12}$, $n\equiv 4\pmod {12}$;
\item[d)] $u\equiv 0,3,4,$ or $7\pmod {12}$, if $n\equiv 9\pmod {12}$.
\end{itemize}


\begin{prop} \label{60}
For any $n=60k+5r$, $r=0,5,8,9$, there exists a decomposition of $K_{n+u}\setminus K_n$ into $3$-suns for every admissible  $u\geq 24k+2r+1$.
\end{prop}

{\em Proof.} Let  $X=\{\infty _1,\infty _2, \ldots,$ $ \infty _{60k+5r}\}$, $r=0,5,8,9$,  
and $u=24k+2r+1+h$, with $h\geq 0$. Set $h=12s+l$, $0\leq l\leq11$ ($l$ depends on $r$), and distinguish the following cases.

\noindent   \emph{Case $1$}:  $r=0,5,8,9$ and $l=0$ (odd $u$).  \\
Write  $\langle Z_u\cup X, D_u \rangle$=$\langle Z_u\cup X, D \rangle  \cup  \langle Z_u, L \rangle$, where $D= [6s+1, 12k+r+6s]$, $|D|=12k+r$,  and   $L= [1,6s]$, and apply Lemmas  \ref{5-1} and \ref{55}.\\
\emph{Case $2$}:  $r=0,9$ and $l=8$ (odd $u$).  \\
Write  $\langle Z_u\cup X, D_u \rangle$=$
\langle Z_u\cup  \{\infty _1,\infty _2,\infty _3\}, \{2,6s+3,6s+5\} \rangle \cup  
\langle Z_u\cup  \{\infty _4\}, \{1\} \rangle \cup  
\langle Z_u\cup  \{\infty _5\}, \{6s+4\} \rangle \cup  
\langle Z_u\cup (X\setminus \{\infty _1,\infty _2,  \infty _3,\infty _4,\infty _5\}), D' \rangle  \cup  \langle Z_u, L \rangle$, where $D'= [ 6s+6,12k+r+6s+4]$, $|D'|=12k+r-1$,   and   $L=[3,6s+2]$, and apply Lemmas \ref{3-3}, \ref{1-1}, \ref{5-1} and \ref{55}. \\
\emph{Case $3$}:  $r=5,8$ and $l=4$ (odd $u$).  \\
Write  $\langle Z_u\cup X, D_u \rangle$=$
\langle Z_u\cup  \{\infty _1,\infty _2,  \infty _3,\infty _4\}, \{2,4\} \rangle \cup  
\langle Z_u\cup \{\infty _5\}, \{1\} \rangle \cup 
\langle Z_u\cup (X\setminus \{\infty _1,\infty _2,  \infty _3,\infty _4,\infty _5\}), D' \rangle  \cup  \langle Z_u, L \rangle$, where $D'=  [ 6s+3,12k+r+6s+2]\setminus \{6s+4\}$, $|D'|=12k+r-1$,   and   $L= [3,6s+4]\setminus \{4,6s+3\}$, and apply Lemmas \ref{4-2}, \ref{1-1}, \ref{5-1} and \ref{55}. \\
  \emph{Case $4$}:  $r=0,8$ and $l=3$ (even $u$).  \\ 
Write  $\langle Z_u\cup X, D_u \rangle$=$
\langle Z_u\cup  \{\infty _1,\infty _2,  \infty _3\}, \{1,\frac{u}{2}\} \rangle \cup  
\langle Z_u\cup \{\infty _4,\infty _5\}, \{2\} \rangle \cup 
\langle Z_u\cup (X\setminus \{\infty _1,\infty _2,\infty _3,\infty _4 ,\infty _{5}\}), D' \rangle  \cup  \langle Z_u, L \rangle$, where $D'= [6s+3, 12k+r+6s+1]$, $|D'|=12k+r-1$,  and   $L=[3,6s+2]$, and apply Lemmas \ref{3-2u2}, \ref{2-1}, \ref{5-1} and \ref{55}. \\
\emph{Case $5$}:  $r=0$ and $l=11$ (even $u$).  \\
Write  $\langle Z_u\cup X, D_u \rangle$=$
\langle Z_u\cup  \{\infty _1,\infty _2,  \infty _3\}, \{1,2,\frac{u}{2}\} \rangle \cup  
\langle Z_u\cup \{\infty _4,\infty _5\}, \{4,6s+3,6s+5,6s+7\}  \rangle \cup 
\langle Z_u\cup (X\setminus \{\infty _1,\infty _2, \infty _3 , \infty _4, \infty _{5}\}), D' \rangle  \cup  \langle Z_u, L \rangle$, where $D'=[ 6s+6,12k+6s+5]\setminus  \{ 6s+7\}$, $|D'|=12k-1$,  and   $L= [3,6s+4]\setminus \{4,6s+3\}$, and apply Lemmas \ref{3-3u2}, \ref{2-4}, \ref{5-1} and \ref{55}. \\
  \emph{Case $6$}:  $r=5$ and $l=1$ (even $u$).  \\
Write  $\langle Z_u\cup X, D_u \rangle$=$
\langle Z_u\cup  \{\infty _1,\infty _2,  \ldots,\infty _6\}, \{1,\frac{u}{2}\} \rangle \cup  
\langle Z_u\cup \{\infty _7,\infty _8,  \infty _9,$ $\infty _{10}\}, \{2\} \rangle \cup 
\langle Z_u\cup (X\setminus \{\infty _1,\infty _2,  \ldots,\infty _{10}\}), D' \rangle  \cup  \langle Z_u, L \rangle$, where $D'= [6s+3, 12k+6s+5]$, $|D'|=12k+3$,   and   $L=[3,6s+2]$, and apply Lemmas \ref{6-2u2}, \ref{4-1}, \ref{5-1} and \ref{55}. \\
\emph{Case $7$}:  $r=5,9$ and $l=9$ (even $u$).  \\ 
Write  $\langle Z_u\cup X, D_u \rangle$=$
\langle Z_u\cup  \{\infty _1,\infty _2,  \infty _3\}, \{1,\frac{u}{2}\} \rangle \cup  
\langle Z_u\cup \{\infty _4,\infty _5\}, \{2,6s+3,6s+4,6s+5\} \rangle \cup 
\langle Z_u\cup (X\setminus \{\infty _1,\infty _2,\infty _3,\infty _4, \infty _{5}\}), D' \rangle  \cup  \langle Z_u, L \rangle$, where $D'= [6s+6, 12k+r+6s+4]$, $|D'|=12k+r-1$,   and   $L=[3,6s+2]$, and apply Lemmas  \ref{3-2u2}, \ref{2-4}, \ref{5-1} and \ref{55}.  \\
 \emph{Case $8$}:  $r=8$ and $l=7$ (even $u$).  \\
Write  $\langle Z_u\cup X, D_u \rangle$=$
\langle Z_u\cup  \{\infty _1,\infty _2,  \infty _3\}, \{1,2,\frac{u}{2}\} \rangle \cup  
\langle Z_u\cup \{\infty _4\}, \{4\} \rangle \cup  
\langle Z_u\cup \{\infty _5\} , \{6s+5\} \rangle \cup 
\langle Z_u\cup (X\setminus \{\infty _1,\infty _2, \infty _3 , \infty _4, \infty _{5}\}), D' \rangle  \cup  \langle Z_u, L \rangle$, where $D'=  [6s+3, 12k+6s+11]\setminus \{6s+4, 6s+5\}$, $|D'|=12k+7$,   and   $L= [3,6s+4]\setminus \{4,6s+3\}$, and apply Lemmas  \ref{3-3u2}, \ref{1-1}, \ref{5-1} and \ref{55}. \\
\emph{Case $9$}:  $r=9$ and $l=5$ (even $u$).  \\
Write  $\langle Z_u\cup X, D_u \rangle$=$
\langle Z_u\cup  \{\infty _1,\infty _2,  \infty _3\}, \{1,\frac{u}{2}\} \rangle \cup  
\langle Z_u\cup \{\infty _4\}, \{2\} \rangle \cup  
\langle Z_u\cup \{\infty _5\} , \{4\} \rangle \cup 
\langle Z_u\cup (X\setminus \{\infty _1,\infty _2, \infty _3 , \infty _4, \infty _{5}\}), D' \rangle  \cup  \langle Z_u, L \rangle$, where $D'= [6s+3,12k+6s+11]\setminus \{6s+4\}$, $|D'|=12k+8$,  and   $L= [3,6s+4]\setminus \{4,6s+3\}$, and apply Lemmas \ref{3-2u2},  \ref{1-1}, \ref{5-1} and \ref{55}. \hfill$\Box$
\bigskip

\begin{prop} \label{70}
For any $n=60k+5r+1$, $r=0,3,4,7$, there exists a decomposition of $K_{n+u}\setminus K_n$ into $3$-suns for every admissible  $u\geq 24k+2r+2$.
\end{prop}

{\em Proof.} Let  $X=\{\infty _1,\infty _2, \ldots,$ $ \infty _{60k+5r+1}\}$, $r=0,3,4,7$,  
and $u=24k+2r+2+h$, with $h\geq 0$. Set $h=12s+l$, $0\leq l\leq11$, and distinguish the following cases.

\noindent   
\emph{Case $1$}:  $r=0,3$ and $l=1$ (odd $u$). \\
Write  $\langle Z_u\cup X, D_u \rangle$=$
\langle Z_u\cup \{\infty \}, \{6s+2\} \rangle \cup 
\langle Z_u\cup (X\setminus \{\infty \}), D' \rangle  \cup  \langle Z_u, L \rangle$, where $D'= [6s+1,12k+r+6s+1]\setminus \{6s+2\}$, $|D'|=12k+r$,   and   $L= [1,6s]$, and apply Lemmas \ref{1-1}, \ref{5-1} and \ref{55}.\\
\emph{Case $2$}:  $r=0,3,4,7$ and $l=9$ (odd $u$).  \\
Write  $\langle Z_u\cup X, D_u \rangle$=$
\langle Z_u\cup  \{\infty _1,\infty _2,  \infty _3\}, \{1,6s+3,6s+4\} \rangle \cup  
\langle Z_u\cup \{\infty _4,\infty _5,$ $\infty _6\}, \{2,6s+5,6s+7\} \rangle \cup 
\langle Z_u\cup (X\setminus \{\infty _1,\infty _2,  \ldots,\infty _6\}), D' \rangle  \cup  \langle Z_u, L \rangle$, where $D'=[6s+6,12k+r+6s+5]\setminus \{6s+7\}$, $|D'|=12k+r-1$,    and   $L=[3,6s+2]$, and apply Lemmas \ref{3-3}, \ref{5-1}  and \ref{55}.\\
\emph{Case $3$}:  $r=4^{\star},7$ and $l=5$ (odd $u$).  \\
Write  $\langle Z_u\cup X, D_u \rangle$=$
\langle Z_u\cup  \{\infty _1,\infty _2,  \infty _3,\infty _4\}, \{2,4\} \rangle \cup  
\langle Z_u\cup \{\infty _5\}, \{1\} \rangle \cup  
\langle Z_u\cup \{\infty _6\}, \{6s+8\} \rangle\cup 
\langle Z_u\cup (X\setminus \{\infty _1,\infty _2,  \ldots,\infty _6\}), D' \rangle  \cup  \langle Z_u, L \rangle$, where $D'= [6s+3,12k+r+6s+3]\setminus  \{6s+4,6s+8\}$, $|D'|=12k+r-1$,   and   $L= [3,6s+4]\setminus \{4,6s+3\}$, and apply Lemmas \ref{4-2}, \ref{1-1}, \ref{5-1} and \ref{55}.\\
\emph{Case $4$}:  $r=0,4$ and $l=6$ (even $u$).  \\
Write  $\langle Z_u\cup X, D_u \rangle$=$
\langle Z_u\cup  \{\infty _1,\infty _2,  \infty _3\}, \{1,\frac{u}{2}\} \rangle \cup  
\langle Z_u\cup \{\infty _4,\infty _5,\infty _6\}, \{2,6s+3, 6s+5\} \rangle \cup 
\langle Z_u\cup (X\setminus \{\infty _1,\infty _2,\ldots,\infty _6\}), D' \rangle  \cup  \langle Z_u, L \rangle$, where $D'= [6s+4,12k+r+6s+3]\setminus  \{6s+5\}$, $|D'|=12k+r-1$,    and   $L=[3,6s+2]$, and apply Lemmas \ref{3-2u2}, \ref{3-3}, \ref{5-1} and \ref{55}.  \\
\emph{Case $5$}:  $r=0$ and $l=10$ (even $u$).  \\
Write  $\langle Z_u\cup X, D_u \rangle$=$
\langle Z_u\cup  \{\infty _1,\infty _2,  \ldots,\infty _6\}, \{1,\frac{u}{2}\} \rangle \cup  
\langle Z_u\cup \{\infty _7,\infty _8,  \infty _9,$ $\infty _{10}\}, \{2\} \rangle \cup  
\langle Z_u\cup \{\infty _{11}\}, \{4, 6s+3,6s+5,6s+6, 6s+7\} \rangle\cup 
\langle Z_u\cup (X\setminus \{\infty _1,\infty _2,  \ldots,\infty _{11}\}), D' \rangle  \cup  \langle Z_u, L \rangle$, where $D'= [6s+8, 12k+6s+5]$, $|D'|=12k-2$,    and   $L= [3,6s+4]\setminus \{4,6s+3\}$, and apply Lemmas \ref{6-2u2}, \ref{4-1},  \ref{1-5}, \ref{5-1}  and \ref{55}. \\
\emph{Case $6$}:  $r=3,7$ and $l=0$ (even $u$).  \\
Write  $\langle Z_u\cup X, D_u \rangle$=$
\langle Z_u\cup  \{\infty _1,\infty _2, \ldots, \infty _6\}, \{1,\frac{u}{2}\} \rangle  \cup 
\langle Z_u\cup (X\setminus \{\infty _1,\infty _2,\ldots,$ $\infty _6\}), D' \rangle  \cup  \langle Z_u, L \rangle$, where $D'= \{2\}\cup [6s+3,12k+r+6s]$, $|D'|=12k+r-1$,    and   $L=[3,6s+2]$, and apply Lemmas \ref{6-2u2},  \ref{5-1} and \ref{55}.  \\
\emph{Case $7$}:  $r=3$ and $l=4$ (even $u$).  \\
Write  $\langle Z_u\cup X, D_u \rangle$=$
\langle Z_u\cup  \{\infty _1,\infty _2,  \infty _3,  \infty _4\}, \{1,\frac{u}{2}\} \rangle \cup  
\langle Z_u\cup \{\infty _5\} , \{2\} \rangle \cup 
\langle Z_u\cup \{\infty _6\}, \{6s+5\} \rangle \cup  
\langle Z_u\cup (X\setminus \{\infty _1,\infty _2,  \ldots, \infty _{6}\}), D' \rangle  \cup  \langle Z_u, L \rangle$, where $D'= [6s+3, 12k+6s+5]\setminus \{6s+5\}$,  $|D'|=12k+2$,    and   $L=[3,6s+2]$, and apply Lemmas \ref{4-2u2},  \ref{1-1}, \ref{5-1}  and \ref{55}. \\
\emph{Case $8$}:  $r=4$ and $l=2$ (even $u$).  \\
Write  $\langle Z_u\cup X, D_u \rangle$=$
\langle Z_u\cup  \{\infty _1,\infty _2,  \infty _3,  \infty _4\}, \{1,\frac{u}{2}\} \rangle \cup  
\langle Z_u\cup \{\infty _5,\infty _6\}, \{2\} \rangle \cup 
\langle Z_u\cup (X\setminus \{\infty _1,\infty _2,\ldots,\infty _{6}\}), D' \rangle  \cup  \langle Z_u, L \rangle$, where $D'= [6s+3, 12k+6s+5]$, $|D'|=12k+3$,    and   $L=[3,6s+2]$, and apply Lemmas  \ref{4-2u2}, \ref{2-1}, \ref{5-1} and \ref{55}. \\
\emph{Case $9$}:  $r=7$ and $l=8$ (even $u$).  \\
 Write  $\langle Z_u\cup X, D_u \rangle$=$
\langle Z_u\cup  \{\infty _1,\infty _2,  \infty _3\}, \{1,2,\frac{u}{2}\} \rangle \cup  
\langle Z_u\cup \{\infty _4,\infty _5,\infty _6\}, \{4,$ $ 6s+3,6s+7\} \rangle \cup 
\langle Z_u\cup (X\setminus \{\infty _1,\infty _2,\ldots,\infty _{6}\}), D' \rangle  \cup  \langle Z_u, L \rangle$, where $D'= [6s+5,  12k+6s+11]\setminus \{6s+7\}$, $|D'|=12k+6$,    and   $L= [3,6s+4]\setminus \{4,6s+3\}$, and apply Lemmas \ref{3-3u2},  \ref{3-3}, \ref{5-1}  and \ref{55}. 
\hfill$\Box$

\begin{prop} \label{80}
For any $n=60k+5r+2$, $r=2,7,10,11$, there exists a decomposition of $K_{n+u}\setminus K_n$ into $3$-suns for every admissible  $u\geq 24k+2r+2$.
\end{prop}

{\em Proof.} Let  $X=\{\infty _1,\infty _2, \ldots,$ $ \infty _{60k+5r+2}\}$, $r=2,7,10,11$,  
and $u=24k+2r+2+h$, with $h\geq 0$. Set $h=12s+l$, $0\leq l\leq11$,  and distinguish the following cases.

\noindent   
\emph{Case $1$}:  $r=2,11$ and $l=3$ (odd $u$).  \\
Write  $\langle Z_u\cup X, D_u \rangle$=$
\langle Z_u\cup \{\infty_1 \}, \{6s+2\} \rangle \cup \langle Z_u\cup \{\infty_2 \}, \{6s+4\} \rangle \cup 
\langle Z_u\cup (X\setminus \{\infty_1,\infty_2 \}), D' \rangle  \cup  \langle Z_u, L \rangle$, where $D'= [6s+1,  12k+r+6s+2]\setminus  \{6s+2, 6s+4\}$, $|D'|=12k+r$,    and   $L= [1,6s]$, and apply Lemmas  \ref{1-1}, \ref{5-1} and \ref{55}.\\
\emph{Case $2$}:  $r=2,7,10,11$ and $l=7$ (odd $u$).  \\
Write  $\langle Z_u\cup X, D_u \rangle$=$
\langle Z_u\cup  \{\infty _1,\infty _2\}, \{1,2,6s+3,6s+4\} \rangle 
\langle Z_u\cup (X\setminus \{\infty _1,\infty _2\}), $ $ D' \rangle  \cup  \langle Z_u, L \rangle$, where $D'= [6s+5,12k+r+6s+4]$, $|D'|=12k+r$,    and   $L=[3,6s+2]$, and apply Lemmas \ref{2-4}, \ref{5-1} and \ref{55}.\\
\emph{Case $3$}:  $r=7,10$ and $l=11$ (odd $u$).  \\
Write  $\langle Z_u\cup X, D_u \rangle$=$
\langle Z_u\cup  \{\infty _1,\infty _2,  \infty _3\}, \{1,6s+3,6s+4\} \rangle \cup  
\langle Z_u\cup \{\infty _4,\infty _5,$ $\infty _6\}, \{2,6s+5,6s+7\} \rangle \cup  
\langle Z_u\cup \{\infty _7\}, \{6s+8\} \rangle\cup 
\langle Z_u\cup (X\setminus \{\infty _1,\infty _2,  \ldots,\infty _7\}), $ $D' \rangle  \cup  \langle Z_u, L \rangle$, where $D'= [6s+6,12k+r+6s+6]\setminus \{6s+7,6s+8\}$, $|D'|=12k+r-1$,    and   $L=[3,6s+2]$, and apply Lemmas \ref{3-3}, \ref{1-1}, \ref{5-1} and \ref{55}.\\
\emph{Case $4$}:  $r=2$ and $l=6$ (even $u$).  \\
Write  $\langle Z_u\cup X, D_u \rangle$=$
\langle Z_u\cup  \{\infty _1,\infty _2,  \infty _3,  \infty _4\}, \{1,\frac{u}{2}\} \rangle \cup  
\langle Z_u\cup \{\infty _5,\infty _6,\infty _7\}, \{2,$ $ 6s+3,6s+5\} \rangle \cup  
\langle Z_u\cup (X\setminus \{\infty _1,\infty _2,  \ldots, \infty _{7}\}), D' \rangle  \cup  \langle Z_u, L \rangle$, where $D'= [6s+4,12k+6s+5]\setminus  \{6s+5\}$, $|D'|=12k+1$,    and   $L=[3,6s+2]$, and apply Lemmas  \ref{4-2u2}, \ref{3-3},  \ref{5-1} and \ref{55}. \\
\emph{Case $5$}:  $r=2,10$ and $l=10$ (even $u$).  \\
Write  $\langle Z_u\cup X, D_u \rangle$=$
\langle Z_u\cup  \{\infty _1,\infty _2, \ldots,  \infty _6\}, \{1,\frac{u}{2}\} \rangle \cup  
\langle Z_u\cup \{\infty _7\}, \{2, 6s+3,6s+4,6s+5,6s+6\} \rangle \cup  
\langle Z_u\cup (X\setminus \{\infty _1,\infty _2,  \ldots, \infty _{7}\}), D' \rangle  \cup  \langle Z_u, L \rangle$, where $D'=  [6s+7, 12k+r+6s+5]$, $|D'|=12k+r-1$,    and   $L=[3,6s+2]$, and apply Lemmas  \ref{6-2u2}, \ref{1-5}, \ref{5-1}  and \ref{55}. \\
\emph{Case $6$}:  $r=7,11$ and $l=4$ (even $u$).  \\
Write  $\langle Z_u\cup X, D_u \rangle$=$
\langle Z_u\cup  \{\infty _1,\infty _2,  \infty _3\}, \{1,\frac{u}{2}\} \rangle \cup  
\langle Z_u\cup \{\infty _4,\infty _5,  \infty _6, \infty _{7}\}, $ $ \{2,4\} \rangle \cup 
\langle Z_u\cup (X\setminus \{\infty _1,\infty _2,  \ldots,\infty _{7}\}), D' \rangle  \cup  \langle Z_u, L \rangle$, where $D'= [6s+3,12k+r+6s+2]\setminus  \{6s+4\}$, $|D'|=12k+r-1$,    and   $L= [3,6s+4]\setminus \{4,6s+3\}$, and apply Lemmas \ref{3-2u2}, \ref{4-2}, \ref{5-1}  and \ref{55}. \\
\emph{Case $7$}:  $r=7$ and $l=8$ (even $u$).  \\
Write  $\langle Z_u\cup X, D_u \rangle$=$
\langle Z_u\cup  \{\infty _1,\infty _2,\infty _3\}, \{1,\frac{u}{2}\} \rangle \cup  \{\infty _4,\infty _5,\infty _6\}, \{2, 6s+3,6s+5\} \rangle
\langle Z_u\cup \{\infty _7\}, \{6s+7\} \rangle \cup  
\langle Z_u\cup (X\setminus \{\infty _1,\infty _2,  \ldots, \infty _{7}\}), D' \rangle  \cup  \langle Z_u, L \rangle$, where $D'=  [6s+4,12k+6s+11]\setminus \{6s+5,6s+7\}$, $|D'|=12k+6$,     and   $L=[3,6s+2]$, and apply Lemmas  \ref{3-2u2}, \ref{3-3},  \ref{1-1},  \ref{5-1} and \ref{55}. \\
\emph{Case $8$}:  $r=10$ and $l=2$ (even $u$).  \\
Write  $\langle Z_u\cup X, D_u \rangle$=$
\langle Z_u\cup  \{\infty _1,\infty _2, \ldots,  \infty _6\}, \{1,\frac{u}{2}\} \rangle \cup  
\langle Z_u\cup \{\infty _7\}, \{2\} \rangle \cup  
\langle Z_u\cup (X\setminus \{\infty _1,\infty _2,  \ldots, \infty _{7}\}), D' \rangle  \cup  \langle Z_u, L \rangle$, where $D'=  [6s+3,12k+6s+11]$, $|D'|=12k+9$,   and   $L=[3,6s+2]$, and apply Lemmas \ref{6-2u2}, \ref{1-1}, \ref{5-1}  and \ref{55}. \\
\emph{Case $9$}:  $r=11$ and $l=0$ (even $u$).  \\
Write  $\langle Z_u\cup X, D_u \rangle$=$
\langle Z_u\cup  \{\infty _1,\infty _2, \ldots,  \infty _7\}, \{1,\frac{u}{2}\} \rangle  \cup  
\langle Z_u\cup (X\setminus \{\infty _1,\infty _2,  $ $ \ldots, \infty _{7}\}), D' \rangle  \cup  \langle Z_u, L \rangle$, where $D'= \{2\}\cup [6s+3,12k+6s+11]$, $|D'|=12k+10$,   and   $L=[3,6s+2]$, and apply Lemmas  \ref{7-2u2},  \ref{5-1} and \ref{55}.   \hfill$\Box$

\begin{prop} \label{90}
For any $n=60k+5r+3$, $r=2,5,6,9$, there exists a decomposition of $K_{n+u}\setminus K_n$ into $3$-suns for every admissible  $u\geq 24k+2r+3$.
\end{prop}

{\em Proof.} Let  $X=\{\infty _1,\infty _2, \ldots,$ $ \infty _{60k+5r+3}\}$, $r=2,5,6,9$,   
and $u=24k+2r+3+h$, with $h\geq 0$. Set $h=12s+l$, $0\leq l\leq11$,  and distinguish the following cases.

\noindent   
\emph{Case $1$}:  $r=2,5,6,9$ and $l=4$ (odd $u$).  \\
Write  $\langle Z_u\cup X, D_u \rangle$=$
\langle Z_u\cup  \{\infty _1,\infty _2, \infty _3\}, \{1,6s+3,6s+4\} \rangle  \cup  
\langle Z_u\cup (X\setminus \{\infty _1,\infty _2,  $ $\infty _{3}\}), D' \rangle  \cup  \langle Z_u, L \rangle$, where $D'= \{2\}\cup  [6s+5,12k+r+6s+3]$, $|D'|=12k+r$,   and   $L=[3,6s+2]$, and apply Lemmas  \ref{3-3},  \ref{5-1} and \ref{55}.  \\
\emph{Case $2$}:  $r=2,5$ and $l=8$ (odd $u$).  \\
Write  $\langle Z_u\cup X, D_u \rangle$=$
\langle Z_u\cup  \{\infty _1,\infty _2\}, \{1,6s+3,6s+4,6s+5\} \rangle \cup  
\langle Z_u\cup \{\infty _3\}, \{2\} \rangle \cup (X\setminus \{\infty _1,\infty _2,  $ $\infty _{3}\}), D' \rangle  \cup  \langle Z_u, L \rangle$, where $D'=  [6s+6, 12k+r+6s+5]$, $|D'|=12k+r$,   and   $L=[3,6s+2]$, and apply Lemmas \ref{2-4},  \ref{1-1},  \ref{5-1}  and \ref{55}.  \\
\emph{Case $3$}:  $r=6,9$ and $l=0$ (odd $u$).  \\
If $s=0$, then write  $\langle Z_u\cup X, D_u \rangle$=$
\langle Z_u\cup  \{\infty _1,\infty _2,\ldots, \infty _8\}, \{1,\frac{u}{3}\} \rangle  \cup  
\langle Z_u\cup (X\setminus \{\infty _1,\infty _2, \ldots, \infty _8\}), D' \rangle $, where $D'=  [2, 12k+r+1]\setminus \{\frac{u}{3}\}$, $|D'|=12k+r-1$, and apply Lemmas \ref{8-2} and \ref{5-1}.
 If $s>0$, then write  $\langle Z_u\cup X, D_u \rangle$=$
\langle Z_u\cup  \{\infty _1,\infty _2, \infty _3\}, \{1,5s,5s+1\} \rangle  \cup  
\langle Z_u\cup  \{\infty _4,\infty _5, \infty _6\}, \{2,6s+1,6s+3\} \rangle  \cup  
\langle Z_u\cup  \{\infty _7\}, \{6s+2\} \rangle  \cup  
\langle Z_u\cup  \{\infty _8\}, \{6s+4\} \rangle  \cup  
\langle Z_u\cup (X\setminus \{\infty _1,\infty _2,  \ldots,\infty _{8}\}), D' \rangle  \cup  \langle Z_u, L \rangle$, where $D'= \{2s+1,4s\}\cup  [6s+5,12k+r+6s+1]$, $|D'|=12k+r-1$,   and   $L= [3,6s]\setminus \{2s+1,4s,5s,5s+1\}$, and apply Lemmas \ref{3-3}, \ref{1-1}  and \ref{5-1}  to decompose the first five subgraphs,   while to decompose the last one apply Lemma \ref{55} $i)$ and delete the orbit $(S_0)$.  \\
\emph{Case $4$}:  $r=2,6$ and $l=1$ (even $u$).  \\
Write  $\langle Z_u\cup X, D_u \rangle$=$
\langle Z_u\cup  \{\infty _1,\infty _2, \infty _3\}, \{1,\frac{u}{2}\} \rangle  \cup  
\langle Z_u\cup (X\setminus \{\infty _1,\infty _2, \infty _{3}\}), $ $ D' \rangle  \cup  \langle Z_u, L \rangle$, where $D'= \{2\}\cup  [6s+3,12k+r+6s+1]$, $|D'|=12k+r$,   and   $L=[3,6s+2]$, and apply Lemmas \ref{3-2u2}, \ref{5-1}   and \ref{55}.  \\
\emph{Case $5$}:  $r=2^\star$ and $l=5$ (even $u$).  \\
Write  $\langle Z_u\cup X, D_u \rangle$=$
\langle Z_u\cup  \{\infty _1,\infty _2, \ldots, \infty _6\}, \{1,\frac{u}{2}\} \rangle \cup  
\langle Z_u\cup \{\infty _7,\infty _8,  \infty _9,$ $ \infty _{10}\}, \{2\} \rangle \cup 
\langle Z_u\cup \{\infty _{11},\infty _{12},\infty _{13}\}, \{4,6s+3,6s+7\} \rangle \cup 
\langle Z_u\cup (X\setminus \{\infty _1,\infty _2, $ $ \ldots,\infty _{13}\}), D' \rangle  \cup  \langle Z_u, L \rangle$, where $D'= [6s+5, 12k+6s+5]\setminus \{6s+7\}$, $|D'|=12k$,   and   $L= [3,6s+4]\setminus \{4,6s+3\}$, and apply Lemmas  \ref{6-2u2}, \ref{4-1}, \ref{3-3}, \ref{5-1}  and \ref{55}. \\
\emph{Case $6$}:  $r=5,9$ and $l=7$ (even $u$).  \\
Write  $\langle Z_u\cup X, D_u \rangle$=$
\langle Z_u\cup  \{\infty _1,\infty _2,\ldots, \infty _6\}, \{1,\frac{u}{2}\} \rangle  \cup  
\langle Z_u\cup  \{\infty _7,\infty _8\}, \{2,6s+3,6s+4, 6s+5\} \rangle  \cup  \langle Z_u\cup (X\setminus \{\infty _1,\infty _2,\ldots, \infty _{8}\}),   D' \rangle  \cup  \langle Z_u, L \rangle$, where $D'= [6s+6,12k+r+6s+4]$, $|D'|=12k+r-1$,   and   $L=[3,6s+2]$, and apply Lemmas \ref{6-2u2}, \ref{2-4},  \ref{5-1}  and \ref{55}.  \\
\emph{Case $7$}:  $r=5$ and $l=11$ (even $u$).  \\
Write  $\langle Z_u\cup X, D_u \rangle$=$
\langle Z_u\cup  \{\infty _1,\infty _2, \infty _3, \infty _4\}, \{1,\frac{u}{2}\} \rangle  \cup  
\langle Z_u\cup  \{\infty _5,\infty _6\}, \{2,6s+3, 6s+5,6s+6\} \rangle  \cup  
\langle Z_u\cup  \{\infty _7\}, \{4\} \rangle
\cup  
\langle Z_u\cup  \{\infty _8\}, \{6s+7\} \rangle
\cup  \langle Z_u\cup (X\setminus \{\infty _1,\infty _2,\ldots, \infty _{8}\}),   D' \rangle  \cup  \langle Z_u, L \rangle$, where $D'=[6s+8,12k+6s+11]$, $|D'|=12k+4$,   and   $L= [3,6s+4]\setminus \{4,6s+3\}$, and apply Lemmas  \ref{4-2u2}, \ref{2-4}, \ref{1-1}, \ref{5-1} and \ref{55}.  \\
\emph{Case $8$}:  $r=6$ and $l=9$ (even $u$).  \\
Write  $\langle Z_u\cup X, D_u \rangle$=$
\langle Z_u\cup  \{\infty _1,\infty _2, \infty _3\}, \{1,\frac{u}{2}\} \rangle  \cup  
\langle Z_u\cup  \{\infty _4,\infty _5,\infty _6\}, \{2,6s+3, 6s+5\} \rangle  \cup  
\langle Z_u\cup  \{\infty _7\}, \{4\} \rangle
\cup  
\langle Z_u\cup  \{\infty _8\}, \{6s+7\} \rangle
\cup  \langle Z_u\cup (X\setminus \{\infty _1,\infty _2,\ldots,$ $ \infty _{8}\}),   D' \rangle  \cup  \langle Z_u, L \rangle$, where $D'=[6s+6,12k+6s+11]\setminus \{ 6s+7\}$, $|D'|=12k+5$,    and   $L= [3,6s+4]\setminus \{4,6s+3\}$, and apply Lemmas \ref{3-2u2}, \ref{3-3}, \ref{1-1}, \ref{5-1}  and \ref{55}.  \\
\emph{Case $9$}:  $r=9$ and $l=3$ (even $u$).  \\
Write  $\langle Z_u\cup X, D_u \rangle$=$
\langle Z_u\cup  \{\infty _1,\infty _2, \infty _3\}, \{1,2,\frac{u}{2}\} \rangle  \cup  
\langle Z_u\cup (X\setminus \{\infty _1,\infty _2, \infty _{3}\}), $ $ D' \rangle  \cup  \langle Z_u, L \rangle$, where $D'= [6s+3,12k+6s+11]$, $|D'|=12k+9$,    and   $L=[3,6s+2]$, and apply Lemmas \ref{3-3u2},  \ref{5-1}  and \ref{55}.   \hfill$\Box$

\begin{prop} \label{100}
For any $n=60k+5r+4$, $r=0,1,4,9$, there exists a decomposition of $K_{n+u}\setminus K_n$ into $3$-suns for every admissible  $u\geq 24k+2r+3$.
\end{prop}

{\em Proof.} Let  $X=\{\infty _1,\infty _2, \ldots,$ $ \infty _{60k+5r+4}\}$, $r=0,1,4,9$,  
and $u=24k+2r+3+h$, with $h\geq 0$. Set $h=12s+l$, $0\leq l\leq11$,  and distinguish the following cases.

\noindent   
\emph{Case $1$}:  $r=0,1^\star,4,9$ and $l=2$ (odd $u$).  \\
Write  $\langle Z_u\cup X, D_u \rangle$=$
\langle Z_u\cup  \{\infty _1,\infty _2,  \infty _3,  \infty _4\}, \{2,4\} \rangle  \cup 
\langle Z_u\cup (X\setminus \{\infty _1,\infty _2,  \infty _3, $ $ \infty _4\}), D' \rangle  \cup  \langle Z_u, L \rangle$, where $D'= \{1,6s+3\}\cup [6s+5,12k+r+6s+2]$, $|D'|=12k+r$,    and   $L= [3,6s+4]\setminus \{4,6s+3\}$, and apply Lemmas \ref{4-2},  \ref{5-1} and \ref{55}. \\
\emph{Case $2$}:  $r=0,9$ and $l=6$ (odd $u$).  \\
Write  $\langle Z_u\cup X, D_u \rangle$=$
\langle Z_u\cup  \{\infty _1,\infty _2, \infty _3\}, \{1,6s+3,6s+4\} \rangle  \cup  \langle Z_u\cup  \{\infty _4\}, \{2\} \rangle  \cup  
\langle Z_u\cup (X\setminus \{\infty _1,\infty _2,  \infty _{3} , \infty _4\}), D' \rangle  \cup  \langle Z_u, L \rangle$, where $D'= [6s+5,12k+r+6s+4]$, $|D'|=12k+r$,    and   $L=[3,6s+2]$, and apply Lemmas  \ref{3-3}, \ref{1-1}, \ref{5-1}  and \ref{55}.  \\
\emph{Case $3$}:  $r=1,4$ and $l=10$ (odd $u$).  \\
Write  $\langle Z_u\cup X, D_u \rangle$=$
\langle Z_u\cup  \{\infty _1,\infty _2\}, \{1,6s+3,6s+5,6s+6\} \rangle  \cup  
\langle Z_u\cup  \{\infty _3\}, \{2\} \rangle  \cup  
\langle Z_u\cup  \{\infty _4\}, \{6s+4\} \rangle
\cup  \langle Z_u\cup (X\setminus \{\infty _1,\infty _2, \infty _{3}, \infty_4\}),   D' \rangle  \cup  \langle Z_u, L \rangle$, where $D'= [6s+7, 12k+r+6s+6]$, $|D'|=12k+r$,    and   $L=[3,6s+2]$, and apply Lemmas \ref{2-4}, \ref{1-1}, \ref{5-1}  and \ref{55}.  \\
\emph{Case $4$}:  $r=0,4$ and $l=5$ (even $u$).  \\
Write  $\langle Z_u\cup X, D_u \rangle$=$
\langle Z_u\cup  \{\infty _1,\infty _2,\ldots, \infty _6\}, \{1,\frac{u}{2}\} \rangle  \cup  \langle Z_u\cup  \{\infty _7,\infty _8, \infty _9\}, \{2, $ $6s+3,6s+5\} \rangle  \cup  
\langle Z_u\cup (X\setminus \{\infty _1,\infty _2, \ldots, \infty_{9}\}),  D' \rangle  \cup  \langle Z_u, L \rangle$, where $D'=  [6s+4,12k+r+6s+3]\setminus\{6s+5\}$, $|D'|=12k+r-1$,    and   $L=[3,6s+2]$, and apply Lemmas \ref{6-2u2}, \ref{3-3}, \ref{5-1}    and \ref{55}.  \\
\emph{Case $5$}:  $r=0$ and $l=9$ (even $u$).  \\
Write  $\langle Z_u\cup X, D_u \rangle$=$
\langle Z_u\cup  \{\infty _1,\infty _2, \infty _3, \infty _4\}, \{1,\frac{u}{2}\} \rangle  \cup  
\langle Z_u\cup  \{\infty _5,\infty _6,\infty _7\}, $ $\{2,6s+3, 6s+5\} \rangle  \cup  
\langle Z_u\cup  \{\infty _8\}, \{4\} \rangle
\cup  
\langle Z_u\cup  \{\infty _9\}, \{6s+7\} \rangle
\cup  \langle Z_u\cup (X\setminus \{\infty _1,\infty _2,\ldots, \infty _{9}\}),   D' \rangle  \cup  \langle Z_u, L \rangle$, where $D'= [6s+6,12k+6s+5]\setminus\{6s+7\}$, $|D'|=12k-1$,    and   $L=[3,6s+4]\setminus \{4,6s+3\}$, and apply Lemmas \ref{4-2u2} ,  \ref{3-3}, \ref{1-1}, \ref{5-1} and \ref{55}.  \\
\emph{Case $6$}:  $r=1$ and $l=7$ (even $u$).  \\
Write  $\langle Z_u\cup X, D_u \rangle$=$
\langle Z_u\cup  \{\infty _1,\infty _2, \ldots, \infty _7\}, \{1,\frac{u}{2}\} \rangle  \cup  
\langle Z_u\cup  \{\infty _8,\infty _9\}, $ $\{2,4,$ $6s+3,6s+5\} \rangle  
\cup  \langle Z_u\cup (X\setminus \{\infty _1,\infty _2,$ $\ldots, \infty _{9}\}),   D' \rangle  \cup  \langle Z_u, L \rangle$, where $D'= [ 6s+6,12k+6s+5]$, $|D'|=12k$,    and   $L= [3,6s+4]\setminus \{4,6s+3\}$, and apply Lemmas \ref{7-2u2}, \ref{2-4},  \ref{5-1}  and \ref{55}.  \\
\emph{Case $7$}:  $r=1,9$ and $l=11$ (even $u$).  \\
Write  $\langle Z_u\cup X, D_u \rangle$=$
\langle Z_u\cup  \{\infty _1,\infty _2, \infty _3\}, \{1,\frac{u}{2}\} \rangle   \cup 
\langle Z_u\cup \{\infty _{4}\}, \{2,4,6s+3,6s+5,6s+6\} \rangle 
\cup  \langle Z_u\cup (X\setminus \{\infty _1,\infty _2, \infty _{3}, \infty _{4}\}),   D' \rangle  \cup  \langle Z_u, L \rangle$, where $D'= [ 6s+7,12k+r+6s+6]$, $|D'|=12k+r$,    and   $L= [3,6s+4]\setminus \{4,6s+3\}$, and apply Lemmas \ref{3-2u2}, \ref{1-5},  \ref{5-1}  and \ref{55}.  \\
\emph{Case $8$}:  $r=4$ and $l=1$ (even $u$).  \\
Write  $\langle Z_u\cup X, D_u \rangle$=$
\langle Z_u\cup  \{\infty _1,\infty _2, \infty _3, \infty _4\}, \{1,\frac{u}{2}\} \rangle  \cup  
\langle Z_u\cup (X\setminus \{\infty _1,\infty _2, $ $\infty _{3}, \infty _4\}),   D' \rangle  \cup  \langle Z_u, L \rangle$, where $D'= \{2\}\cup [6s+3,12k+6s+5]$, $|D'|=12k+4$,    and   $L=[3,6s+2]$, and apply Lemmas \ref{4-2u2},  \ref{5-1}  and \ref{55}.  \\
\emph{Case $9$}:  $r=9$ and $l=3$ (even $u$).  \\
Write  $\langle Z_u\cup X, D_u \rangle$=$
\langle Z_u\cup  \{\infty _1,\infty _2, \infty _3\}, \{1,\frac{u}{2}\} \rangle  \cup  \langle Z_u\cup  \{\infty _4\}, \{2\} \rangle  \cup  
\langle Z_u\cup (X\setminus \{\infty _1,\infty _2, $ $\infty _{3}, \infty _4\}),   D' \rangle  \cup  \langle Z_u, L \rangle$, where $D'= [6s+3,12k+6s+11]$, $|D'|=12k+9$,   and   $L= [3,6s+2]$, and apply Lemmas \ref{3-2u2}, \ref{1-1},   \ref{5-1} and \ref{55}.   \hfill$\Box$

\bigskip 

Combining Lemma \ref{NC} and Propositions \ref{60}--\ref{100}
gives our main theorem.

\begin{theo}  
Any 3SS$(n)$ can be embedded in a 3SS$(m)$ if and only if $m\geq
\frac 75 n+1$ or $m=n$.
\end{theo}

\noindent {\Large \textbf{Appendix}}

\begin{itemize}
 \item $n=21$, $u=12s+15$\\
Write  $\langle Z_u\cup X, D_u \rangle$=$
\langle Z_u\cup  \{\infty _1,\infty _2,  \infty _3,\infty _4\}, \{2,4\} \rangle \cup  
\langle Z_u\cup \{\infty _5\}, \{1\} \rangle \cup  
\langle Z_u\cup \{\infty _6\}, \{6s+7\} \rangle\cup 
\langle Z_u\cup (X\setminus \{\infty _1,\infty _2,  \ldots,\infty _6\}), \{6s+3,6s+5,6s+6\} \rangle  \cup  \langle Z_u, L \rangle$, where  $L= [3,6s+4]\setminus \{4,6s+3\}$, and apply Lemmas \ref{4-2}, \ref{1-1}, \ref{5-1} and \ref{55}.
\item $n=13$, $u=12s+12$ \\
Write  $\langle Z_u\cup X, D_u \rangle$=$
\langle Z_u\cup  \{\infty _1,\infty _2, \ldots, \infty _6\}, \{1,6s+6\} \rangle \cup  
\langle Z_u\cup \{\infty _7,$ $ \infty _8,  \infty _9,\infty _{10}\}, \{2\} \rangle \cup 
\langle Z_u\cup \{\infty _{11},\infty _{12},\infty _{13}\}, \{4,6s+3,6s+5\} \rangle \cup  \langle Z_u, L \rangle$, where   $L= [3,6s+4]\setminus \{4,6s+3\}$, and apply Lemmas  \ref{6-2u2}, \ref{4-1}, \ref{3-3}  and \ref{55}.
\item $n=9$, $u=12s+7$ \\
Write  $\langle Z_u\cup X, D_u \rangle$=$
\langle Z_u\cup  \{\infty _1,$ $\infty _2,  \infty _3,  \infty _4\}, \{2,4\} \rangle  \cup 
\langle Z_u\cup \{\infty _5,\infty _6, $ $ \infty _7, \infty _8, \infty _9\}, \{1\}\rangle  \cup  \langle Z_u, L \rangle$, where  $L= [3,6s+3]\setminus \{4\}$, and apply Lemmas \ref{4-2},  \ref{5-1} and decompose $ \langle Z_u, L \rangle$ as in Lemma \ref{55} $iii)$, taking in account that $|6s+4|_{12s+7}=6s+3$.
 \end{itemize}

\end{document}